\documentclass[a4paper,10pt]{article}
\def \Z {{\mathbf {Z}}}

\def \N {{\mathbf {N}}}

\def \B {{\cal B}}

\title{ Компактные семейства и типичные  энтропийные инварианты сохраняющих меру действий}
\author{В.В.Рыжиков}
\date{12.02.2021}
\usepackage[T1]{fontenc} 
\usepackage[cp1251]{inputenc}  
\usepackage[russian]{babel}
\begin{document}
\large
\maketitle
\begin{abstract}  
Статья предназначена для выпуска   Трудов Mосковского Mатематического Общества (82:1), \\
посвященного 80-летним юбилеям В.И. Оселедца и А.М. Степина  и  памяти А.М. Степина, \\
 покинувшего  нас в ноябре 2020г.

Содержание. Для компактного множества  действий  энтропия типа энтропии Кушниренко подбирается  
таким образом, чтобы на этом множестве
она равнялась  нулю, но принимала бесконечное значение для типичных действий. 
Как следствие получен результат о том, 
 что типичные сохраняющие меру преобразования не изоморфны   изометрическим 
перекладываниям конечного  набора геометрических фигур.

Библиография: 19 названий, УДК: 517.987, \ MSC: Primary 28D05; Secondary 58F11

Ключевые слова и фразы: \it типичные эргодические действия, 
энтропия Кушниренко, 

компактые семейства динамических систем.\rm

\end{abstract} 

\section{Введение}  Настоящая работа дает отрицательный ответ на  вопрос: \it будут ли типичны преобразования,
изоморфные изометрическому перекладыванию  конечного числа  геометрических фигур? \rm
Задача возникла в связи с результатом  Девис и Чайки \cite{CD}  о нетипичности 
перекладываний конечного числа отрезков.  Перекладывания (interval exchange transformations, IETs) 
давно являются чрезвычайно популярными объектами в эргодической теории. Вспоминая их самые первые  применения, 
 отметим, что в работе Степина \cite{S}  перекладываниями
было реализовано действие диадической группы с  отсутствием группового свойства спектра,
 а   Оселедцем в \cite{O} получен  ряд  результатов о перекладываниях, 
включая реализацию спектров конечной кратности.

В доказательстве    нетипичности перекладываний использовались 
  отклонения от перемешивания, обнаруженные Катком \cite{Ka}. Авторы  \cite{CD}  
сравнили  регулярность  таких отклонений с  нерегулярностью отклонений, характерной 
 для типичных преобразований. 
Отсутствие перемешивания для перекладываний прямоугольников, насколько нам известно, не установлено. 
  Как   доказать их нетипичность? Предлагается метод, использующий  
   небольшую модификацию энтропии Кушниренко \cite{Ku}: 
  перекладывания имеют нулевую специально выбранную  энтропию, а типичные преобразования --  бесконечную энтропию. 
Аналогичные результаты справедливы для изометрических перекладываний конечного набора  геометрических фигур. 
Доказательство их нетипичности   использует равенство нулю классической энтропии (\cite{Bu}).  
Подмножество в группе  сохраняющих меру
преобразований $Aut$, состоящее в точности из классов сопряженности элементов  множества 
$K\subset Aut$, обозначим через $K^{Aut}.$
Вопрос, оставшийся нерешенным,  звучит так: 

\it  верно ли, что множество    $K^{Aut}$ имеет первую категорию в группе ${Aut}$, оснащенной  метрикой Халмоша, 
если $K\subset Aut$ -- компакт? \rm  

 В случае, когда  все  элементы компакта $K$ имеют нулевую классическую энтропию, то ответ, 
как будет показано в заметке,  положительный.  
Это дает большое многообразие нетипичных классов, выходящее  за рамки  перекладываний конечного множества фигур.
Аналогичные  соображения можно с успехом применить к потокам, к групповым действиям и даже 
к  перемешивающим действиям  с метрикой Альперна-Тихонова \cite{ST},  но мы  ограничимся в заметке небольшим кругом задач. 
Любопытно, что преобразования с положительной энтропией, являясь нетипичными,  
 явно "мешают" ответить на поставленный выше вопрос.
Причина, вероятно,  кроется в том, что типичная динамика на чрезвычайно больших интервалах времени   похожа 
на нетипичную (недетерминированную) динамику с малой  положительной энтропией 
(хотя  на других очень больших интервалах времени она  не похожа на нее).

О современных аспектах теории типичных преобразований, в том числе о революционном результате 
 Кинга (King) о корнях,  тонком результате Лазаро (Lazaro) и де ля Рю (de la Rue) 
 о включении типичного преобразования в поток, неожиданных теоремах Агеева 
о факторах и нетривиальных самоприсоединениях(self-joinings) 
см., например,   работы \cite{GTW}, \cite{SE}, \cite{ST}.

\section{ Несколько  типичных свойств преобразований}
 Фиксируем стандартное вероятностное пространство $(X,\B,\mu)$ и рассмотрим 
группу его автоморфизмов $ Aut $,   снабженную полной метрикой Халмоша $\rho $: 
$$ \rho(S,T)=\sum_i \mu(SA_i\Delta TA_i)+\mu(S^{-1}A_i\Delta T^{-1}A_i),$$
где  семейство  множеств  $\{A_i\}$ плотно в алгебре $\B$.
Говорят, что свойство типично, если  множество автоморфизмов (далее мы их называем преобразованиями),
  обладающее этим свойством, 
  содержит  некоторое  $G_\delta$-множество, плотное в  $Aut$.

\vspace{3mm}
\bf 2.1.  Типичность и нетипичность слабых пределов. \rm 
Отсутствие перемешивания и наличие слабого  перемешивания, как показали Рохлин и Халмош,
 являются типичными свойствами. Напомним, как доказывается полезный более общий факт. 

 Функция от оператора  $Q(T)$  называется допустимой,
 если она имеет вид  $$Q(T)=a\Theta+\sum_i a_iT^i, \ \ a, a_i\geq 0, \ \sum_i a_i=1-a,$$  
где $\Theta$ -- oператор ортопроекции на пространство констант в  $ L_2(X,\B, \mu)$.
Если $$T^m\to_w \kappa \Theta + (1-\kappa)I,$$
такое преобразование $T$ называется $\kappa$-перемешивающим. 
Перемешивающим, напомним,  называется преобразование $T$, удовлетворяющее условию
$T^m\to\Theta$ при $m\to\infty$. Читатель, вероятно, уже заметил, что мы одинаково обозначаем преобразование и 
индуцированный преобразованием оператор в $ L_2(X,\B, \mu)$.

Свойство $\kappa$-перемешивания,   сыгравше  важную роль в решении 
задачи Колмогорова о групповом свойстве спектра, 
является типичным свойством, как показал Степин \cite{St}. Определение $\kappa$-перемешивания 
появилось в работе \cite{O1}, где Оселедец высказал
  гипотезу о том, что этим свойством могут обладать перекладывания отрезков 
(гипотеза  через 30 лет была доказана автором).  
В работе \cite{KS}  Каток и Степин рассматривали преобразования 
с $(n,n+1)$-аппроксимацией. Оказывается, такие преобразования
обладают всеми допустимыми слабыми пределами \cite{R1}. 
Большинство перекладываний 3-х отрезков обладает  $(n,n+1)$-аппроксимацией  (см. \cite{CE}).

\vspace{3mm}
\bf Теорема 2.1. \it Для бесконечного множества  $M$ и допустимых функций    $Q$,$R$ 

(i)  множество преобразований $T$  таких , что для некоторого бесконечного подмножества $M'\subset M$
выполнено  $T^m\to_w R(T)$  при     $m\in M'$ , $m\to\infty$, является типичным множеством,

(ii)  множество   преобразований $T$    таких, что  $ T^m\to_w Q(T)$  при     $m\in M, m\to\infty$, 
является множеством первой категории.

\vspace{3mm}
\rm  

В частных случаях при  $R=I$, $Q=\Theta$, получаем упомянутые результаты  Халмоша и Рохлина,
 а подстановка  $R=\kappa\Theta +(1-\kappa)I$  дает  типичность свойства $\kappa$-перемешивания.\rm

\vspace{3mm}
Доказательство. (i).  Фиксируем плотное в $Aut$  множество преобразований   $\{J_q\}$.  Найдется
 бесконечное подмножество  $M'\subset M$ и слабо перемешивающее преобразование $S$ такое, что  $S^m\to_w R(S)$, 
 $m\in M'$.  
Нужное  $S$ несложно реализовать при помощи конструкции  ранга 1 (\cite{R3}).  
Надстройки (spacers) этой конструкции обеспечат 
требуемые  слабые пределы, причем  высоты башен $h_j$ легко заставить принимать значения из $M$. 
Тем самым мы получим   $S^{h_j}\to_w R(S)$,  $h_j\in M$.

Пусть $w$ обозначает метрику, задающую слабую операторную топологию.
 Для любых  $n$ и $q$  найдется  $m=m(n,q)$ и окрестность 
$U(n,q)$ преобразования   $J_q^{-1}SJ_q$ такая, что неравенство 
$$w(T^m,R(T))< \frac 1 n$$ 
выполняется для всех   $T\in U(n,q)$. Получаем  $G_\delta$-множество
$$W=\bigcap_n\bigcup_q U(n,q),$$  плотное в $Aut$, так как класс сопряженности
эргодического преобразования $S$  плотен в $Aut$ (следствие леммы Рохлина-Халмоша). $W$ состоит из преобразовний,
удовлетворяющих условию пункта (i).

 (ii).  Это утверждение логически вытекает из (i). Действительно,  фиксируем бесконечное  $M$,
для типичного преобразования $S$   найдется бесконечное подмножество в $M$, на элементах  которого   
$\{S: S^m\to_w R(S)\neq Q(S)$.
Таким образом,  условие $\{T: T^m\to_w Q(T)$  при     $m\in M, m\to\infty$  выполняется только для нетипичных $T$.

\vspace{3mm}
\bf 2.2. Асимметрия и кратные слабые пределы.  \rm Если $T$ и  $T^{-1}$ не сопряжены в $Aut$, такое преобразование 
$T$ называем  асимметричным.
Типичные преобразования  асимметричны (дель Джунко, \cite{Ju}).  
 В работе \cite{O2} Оселедец построил $\Z_5$-расширение над
перекладыванием отрезков, не изоморфное обратному. Вероятно, аналогичным способом можно получить 
асимметричное перекладывание.
  В  \cite{R2} даны примеры преобразований с очевидной асимметрией, вытекающей из приведенных ниже  
 асимптотических свойств. Заметим, что такими свойствами не могут обладать перекладываний 3-х отрезков,
так как  они являются композицией двух инволюций, следовательно,  они симметричны.  Любопытно узнать, 
бывают ли, например,  несимметричные перекладывания отрезков $I_1, I_2,I_3,I_4$  в порядке $I_3, I_4,I_1,I_2$?
Отметим, что перекладывание в порядке  $I_4, I_3,I_2,I_1$ изоморфно обратному как  
композиция симметрии всего отрезка c инволюцией, образованной одновременными   симметриями  внутри 
каждого из четырех отрезков  $I_4, I_3,I_2,I_1$.

\vspace{3mm}
\bf  Теорема 2.2. \it  Пусть  преобразование $T$ обладает следующим свойством:  для некоторых последовательностей  
$m_i,n_i\to\infty$   для любого измеримого множества $A$ имеют место сходимости
$$\mu(A\cap T^{m_i}A\cap T^{n_i}A)\to (\mu(A)+2\mu(A)^3)/3,$$
$$\mu(A\cap T^{-m_i}A\cap T^{-n_i}A)\to \mu(A)^2.$$
Такие преобразования $T$ асимметричны и типичны.
\rm

\vspace{3mm}
\bf 2.3. Нетипичные неперемешивающие преобразования. \rm

В \cite{R3}  было доказано  следующее утверждение.

\vspace {3mm}
\bf Теорема 2.3. \it Пусть  $ K $ -- компакт в $Aut=Aut(X,\B,\mu)$ и для некоторого  $ r> 0 $ для всех
$ T \in K $ и любого $ n $ найдется  $ m> n $ такое, что  $ w(T^m, \Theta)> r $,
где  $ w $ -- метрика, задающая слабую операторную топологию в $L_2(X,\B,\mu)$.
Тогда множество $K^{Aut}$  имеет первую категорию в $Aut$. \rm

  Если $\Theta$ является единственным  марковским оператором, сплетающим два
преобразования,    то такие преобразования называются дизъюнктными. Это эквивалентно
тому, что они обладают  единственным присоединением (joining) $\mu\times\mu$.   
Дизъюнктные преобразования, неформально говоря,  максимально неизоморфны.  
Обсуждение  предыдущей теоремы  с Бенджи Вейсом (Benjy Weiss) 
 способстовало следующему обобщению результата  Девис и Чайки \cite {CD}.

\vspace{3mm}
\bf Теорема 2.4.  \it Множество преобразований, дизъюнктных со всеми эргодическими перекладываниями
является типичным.\rm

\vspace{3mm}
Доказательство.  Фиксируем плотное в $\B$ семейство множеств $B_i$.
 Перекладывания $p$ отрезков образуют  компакт $K\subset Aut$. В силу частичной жесткости перекладываний,
как известно, найдется положительная константа $c$ (сравнимая с $\frac 1 p$) такая, что 
 для любых  $j$ и  $S\in K$ найдется минимальное число $ N(S,j)>j$ такое, что для некоторого $m$ при 
$j<m\leq N(S,j)$ для всех $i$, $1\leq i\leq j$
выполнено неравенство   
$$\mu(T^mB_i\cap B_i)>c{\mu(B_i)}.\eqno (1)$$
 Так как $K$ -- компакт, $N(S,j)$ будет ограниченной функцией на $K$. 
Обозначим через $N(j)$  максимум значений $N(S,j)$ на $K$.  

Рассмотрим множества $F_j=\{j,j+1,\dots, N(j)\}$.
В \cite{R} доказано, что существует  всюду плотное $G_\delta$-множество  $Y$ такое, что для каждого $T\in Y$
    найдется  последовательность  $F_{j(k)}$, $j(k)\to\infty$ такая, что
$T^{i(k)}\to_w\Theta$ при $i(k)\in F_{j(k)}$, $i(k)\to\infty$.  
 
Из (1) вытекает, что для каждого перекладывания $S$ найдется  последовательность  $\{m_k\}$, 
$m_k\in F_{j(k)}$, для которой  имеет место сходимость
$$   S^{m_k}\to cI+ (1-c) P,$$
где  $P$ -- некоторый марковский оператор, коммутирующий с  $S$.  
Возьмем произвольный $T\in Y$ и покажем, что он дизъюнктен с $S\in K$.  Пусть марковский оператор $J$ сплетает
$S$ и $T$. Имеем
    $$SJ=JT, \ \ S^mJ=JT^m,\ \ S^{m_k}J=JT^{m_k},$$ 
$$cJ+ (1-c) PJ=J\Theta=\Theta.$$
Оператор $\Theta$ является крайней точкой в множестве марковских операторов, сплетающих
эргодический  $S$ со слабо перемешивающим  $T$ (это эквивалентно тому, что $\mu\times\mu$ эргодична относительно 
$S\times T$), поэтому $J=\Theta.$
Это  означает, что  $S$ и  $T$ дизъюнктны.

\section{ Типичность бесконечной энтропии Кушниренко}
Рассмотрим следующую модификацию   энтропии Кушниренко.
 Пусть  $P=\{P_j\}$ -- последовательность конечных подмножеств в счетной бесконечной группе $G$. 
 Для сохраняющего меру действия $T=\{T_g\}$ группы $G$     определим величины 
$$h_j(T,\xi)=\frac 1 {|P_j|}  H\left(\bigvee_{p\in P_j}T_p\xi\right),$$
$$h_{P}(T,\xi)={\limsup_j} \ h_j(T,\xi),$$
$$h_{P}(T)=\sup_\xi h_{P}(T,\xi),$$
где  $\xi$ обозначает конечное измеримое разбиение пространства  $X$, $H(\xi)$ -- энтропия разбиения $\xi$:
$$ H(\{C_1,C_2,\dots, C_n\})=-\sum_{i=1}^n \mu( C_i)\ln \mu( C_i).$$

Нас будет интересовать  только случай  $|P_j|\to\infty$ 
(хотя случай ограниченной мощности также имеет смысл, непосредственно связанный со свойствами типа 
 кратного перемешивания).

\bf 3.1. $P$-энтропия. \rm  Рассмотрим следующий частный случай  $G=Z$, когда  $P_j$ представляют собой
увеличивающиеся в размерах  прогрессии:
  $P_j=\{j,2j,\dots, L(j)j\}$, для некоторой последовательности $L(j)\to\infty$.

\vspace{3mm}
\bf Теорема 3.1. \it Множество $\{S: h_P(S)=\infty\}$ типично.\rm

Доказательство.    Пусть $\{J_q\}$, $q\in \N$, плотно в   $Aut$,  а
$T$ -- бернуллиевское преобразование, обозначим   $T_q=J_q^{-1}TJ_q$.  Множество  $\{T_q\}$ плотно в $Aut$.
 Фиксируем плотное семейство $\{\xi_i\}$ конечных измеримых разбиений.

 Для любых $n,q$ найдется $j=j(n,q)$ такое, что для всех $i\leq n$ выполняется 
$$ h_{j}(T_q,\xi_i)=\frac 1 {L_j} H(\bigvee_{n=1}^{L(j)} T^{nj}_q\xi_i) >H(\xi_i)-\frac 1 n. \eqno (n)$$
В самом деле,  $T_q$ бернуллиевский,  находим разбиение  $\xi$, близкое к фиксированному разбиению $\xi_i$, тогда 
для некоторого числа $m(i,q,n,)$  разбиения  $T^{nj}_q\xi$ независимы при всех $n$ и  $j>m(i,n,q)$.
Это влечет за собой   $(n)$  для всех достаточно больших $j$.

Выбираем окрестность $U(n,q)$ преобразования  $T_q$  такую, что для всех $S\in U(n,q)$ выпоняется неравенство
$$ h_{j}(S,\xi_i) > H(\xi_i)-\frac 1 n.$$ 
Множество 
$$W=\bigcap_n\bigcup_q U(n,q),$$ 
яляется плотным  $G_\delta$. Если  $S\in W$, то для любого $n$ найдется $q(n)$ 
такое, что неравенство    $$ h_{j(n,q(n))}(S,\xi_i)> H(\xi_i)-\frac 1 n$$ выполнено  при $i\leq n$. 
А это приводит к   $h_P(S)=\infty$. Таким образом  множество  $\{S: h_P(S)=\infty\}$ является типичным.

Без особых изменений доказывается следующее более общее утверждение. Автор признателен Льюису Боуэну (Lewis Bowen), 
обратившему  наше внимание на это.

\vspace{3mm}
\bf Теорема 3.1.1. \it Пусть $G$ -- бесконечная счетная группа, а   $P$ -- последовательность множеств
 $P_j=\{p_j(1), p_j(2), \dots, p_j(L(j))\}\subset G$, $L(j)\to\infty$,
  такая, что для любого конечного множества  $F\subset G$ для всех больших  $j$ для всех   $m,n,$  $m\neq n\leq L(j)$,
 произведения  $p_j(m)^{-1}p_j(n)$  не принадлежат $F$ (расстояния между элементами $P_j$ стремятся к $\infty$).

Тогда множество $G$-действий с бесконечной $P$-энтропией  
является типичным множеством в пространстве всех сораняющих меру  $G$-действий.\rm

\vspace{3mm}
{ \bf Компактные семейства с нулевой $P$-энтропией.}
 Обозначим  через  $E_0$ класс преобразований с нулевой энтропией,   
$K^{Aut}=\{J^{-1}SJ:  S\in K, J\in Aut\}$.

\vspace{3mm}
\bf Теорема 3.2. \it Если  $K\subset E_0$ является компактным множеством и  $(Aut,\rho)$, то множество $K^{Aut}$ 
нетипично.\rm

\vspace{3mm}
Доказательство.  Фиксируем плотное семейство конечных разбиений $\xi_i$.
Если    $h(S)=0$, для любого  $i$ имеем
$$h(S^j,\xi)=\lim_{L\to\infty}\frac 1 {L}  H\left(\bigvee_{p=1}^L S^{jp}\xi_i\right)=0.$$
Для  $S\in K$ и  $j$  находим последовательность прогрессий  $P(S)=\{P_{j}(S)\}$
$$   P_j(S)=\{j,2j,\dots, L_S(j)j\}$$  
такую, что 
$$\frac 1 {|P_j(S)|}  H\left(\bigvee_{p\in P_j(S)}S^p\xi_i\right)<\frac 1 j$$
выполняется при  $i<j$. 

 Учитывая структуру множеств $P_j(S)$ и то, что   $K$ -- компакт,  находим последовательность $L(j)\to\infty$ 
такую, что  для любого  $S\in K$ при всех достаточно больших $j$  выполняется $L(j)>L_S(j)$ . 
Тогда для последовательности $P$ из 
расширяющихся прогрессий  $P_j=\{j,2j,\dots, L(j)j\}$ выполнено   $h_P(S)=0$ для всех $S\in K$  и тем самым для всех 
$h_P(T)=0$ для всех $T\in K^{Aut}$. Из теоремы 3.1  вытекает, что  $K^{Aut}$ есть множество первой категории,
что завершает доказательство.

\vspace{3mm}
\bf Вопросы. \rm 
Эргодические преобразования $T$ с положительной энтропией (а также все эргодические  
гауссовские и пуассоновские надстройки) обладают с следующим свойством: найдется апериодическое преобразование $S$ 
и последовательность $n(i)\to\infty$ такие, что имеет место сильная операторная сходимость
$$T^{-n(i)}ST^{n(i)}\to I.$$
Является ли это свойство типичным?

В связи с тем, что типичное преобразование обладает обширной структурой факторов 
(инвариантных сигма-подалгебр),
  возникает  вопрос о значениях $P$-энтропии  для факторов типичного преобразования. 
 Какой из приведенных ниже случаев  типичен:
 $P$-энтропия  нетривиальных факторов преобразования

(i) равна  $\infty$, 

(ii) имеет только значения 0 и $\infty$, 

(iii) 
принимает  все значения от 0 до $\infty$?

\vspace{13mm}
Автор благодарит T.Adams, А. Баштанова, L.Bowen, J.Chaika 
за отклик и полезные замечания и 
особенно J.-P.Thouvenot  и B.Weiss  за плодотворные обсуждения.



E-mail: vryzh@mail.ru
\end{document}